\documentclass[12pt,reqno]{amsart}
\usepackage{amsfonts,amsmath,amssymb}

\allowdisplaybreaks
\advance\textheight by \topskip

\usepackage{mathrsfs}

\usepackage[latin1]{inputenc}
\usepackage{graphicx}

\usepackage{tikz}
\newtheorem{theorem}{Theorem}
\usepackage{euscript}

\def\[{[\! [}
\def\]{]\! ]}

\title[Partial skew Dyck paths]{Partial skew Dyck paths---a kernel method approach}

\author[H.~Prodinger]{Helmut Prodinger}

\address{Helmut Prodinger,
Mathematics Department, Stellenbosch University,
7602 Stellenbosch, South Africa.}
\email{hproding@sun.ac.za}

\date{\today}

\begin{document}

\begin{abstract}
	Skew Dyck are a variation of Dyck paths, where additionally to steps $(1,1)$ and $(1,-1)$ a south-west step $(-1,-1)$ is also allowed, provided that the path
	does not intersect itself. Replacing the south-west step by a red south-east step, we end up with decorated Dyck paths. We analyze partial versions of them
	where the path ends on a fixed level $j$, not necessarily at level 0. We exclusively use generating functions and derive them with the celebrated kernel method.
	
	In the second part of the paper, a dual version is studied, where the paths are read from right to left. In this way, we have two types of up-steps, not two types of
	down-steps, as before.
	
	A last section deals with the variation that the negative territory (below the $x$-axis) is also allowed. Surprisingly, this is more involved in terms of computations.
\end{abstract}

\maketitle

\section{Introduction}

Skew Dyck are a variation of Dyck paths, where additionally to steps $(1,1)$ and $(1,-1)$ a south-west step $(-1,-1)$ is also allowed, provided that the path
does not intersect itself. Otherwise, like for Dyck path, it must never go below the $x$-axis and end eventually (after $2n$ steps) on the $x$-axis.
Here are a few references: \cite{Deutsch-italy, KimStanley, Baril-neu, Prodinger-hex}. The enumerating sequence is
\begin{equation*}
	1, 1, 3, 10, 36, 137, 543, 2219, 9285, 39587, 171369, 751236, 3328218, 14878455,\dots,
\end{equation*}
which is A002212 in \cite{OEIS}.

Skew Dyck appeared very briefly in our recent paper \cite{Prodinger-hex}; here we want to give a more thorough analysis of them, using generating functions and
the kernel method. 
Here is a list of the 10 skew paths consisting of 6 steps:

\begin{figure}[h]
	\begin{equation*}
	\begin{tikzpicture}[scale=0.3]
		\draw [thick](0,0)--(3,3)--(6,0);
	\end{tikzpicture}
	\quad
	\begin{tikzpicture}[scale=0.3]
		\draw [thick](0,0)--(3,3)--(5,1)--(4,0);
	\end{tikzpicture}
	\quad
	\begin{tikzpicture}[scale=0.3]
		\draw [thick](0,0)--(3,3)--(4,2)--(3,1)--(4,0);
	\end{tikzpicture}
	\quad
	\begin{tikzpicture}[scale=0.3]
		\draw [thick](0,0)--(3,3)--(4,2)--(3,1)--(2,0);
	\end{tikzpicture}
	\quad
	\begin{tikzpicture}[scale=0.3]
		\draw [thick](0,0)--(2,2)--(4,0)--(5,1)--(6,0);
	\end{tikzpicture}
\end{equation*}
\begin{equation*}
	\begin{tikzpicture}[scale=0.3]
		\draw [thick](0,0)--(1,1)--(2,0)--(4,2)--(6,0);
	\end{tikzpicture}
	\quad
	\begin{tikzpicture}[scale=0.3]
		\draw [thick](0,0)--(1,1)--(2,0)--(4,2)--(5,1)--(4,0);
	\end{tikzpicture}
	\quad
	\begin{tikzpicture}[scale=0.3]
		\draw [thick](0,0)--(1,1)--(2,2)--(3,1)--(4,2)--(6,0);
	\end{tikzpicture}
	\quad
	\begin{tikzpicture}[scale=0.3]
		\draw [thick](0,0)--(1,1)--(2,2)--(3,1)--(4,2)--(5,1)--(4,0);
	\end{tikzpicture}
	\quad
	\begin{tikzpicture}[scale=0.3]
		\draw [thick](0,0)--(1,1)--(2,0)--(3,1)--(4,0)--(5,1)--(6,0);
	\end{tikzpicture}
\end{equation*}
\caption{All 10 skew Dyck paths of length 6 (consisting of 6 steps).}
\end{figure}
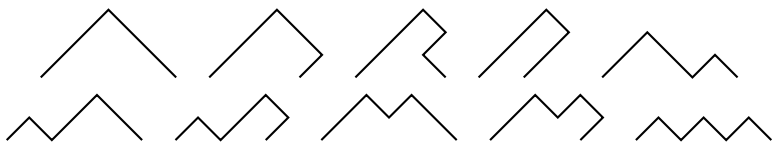

We prefer to work with the equivalent model (resembling more traditional Dyck paths) where
we replace each step $(-1,-1)$ by $(1,-1)$ but label it red. Here is the list of the 10 paths again (Figure 2):

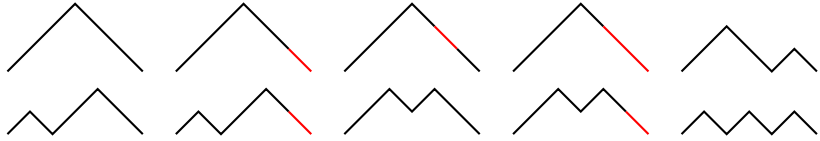
\begin{figure}[h]
\begin{equation*}
	\begin{tikzpicture}[scale=0.3]
		\draw [thick](0,0)--(3,3)--(6,0);
	\end{tikzpicture}
	\quad
	\begin{tikzpicture}[scale=0.3]
		\draw [thick](0,0)--(3,3)--(5,1);
				\draw [thick,red](5,1)--(6,0);
	\end{tikzpicture}
	\quad
	\begin{tikzpicture}[scale=0.3]
		\draw [thick](0,0)--(3,3)--(4,2);
				\draw[red,thick] (4,2)--(5,1);
				\draw [thick](5,1)--(6,0);
	\end{tikzpicture}
	\quad
	\begin{tikzpicture}[scale=0.3]
		\draw [thick](0,0)--(3,3)--(4,2);
		\draw[red,thick](4,2)--(6,0);
	\end{tikzpicture}
	\quad
	\begin{tikzpicture}[scale=0.3]
		\draw [thick](0,0)--(2,2)--(4,0)--(5,1)--(6,0);
	\end{tikzpicture}
\end{equation*}
\begin{equation*}
	\begin{tikzpicture}[scale=0.3]
		\draw [thick](0,0)--(1,1)--(2,0)--(4,2)--(6,0);
	\end{tikzpicture}
	\quad
	\begin{tikzpicture}[scale=0.3]
		\draw [thick](0,0)--(1,1)--(2,0)--(4,2)--(5,1);
		\draw[red,thick] (5,1)--(6,0);
	\end{tikzpicture}
	\quad
	\begin{tikzpicture}[scale=0.3]
		\draw [thick](0,0)--(1,1)--(2,2)--(3,1)--(4,2)--(6,0);
	\end{tikzpicture}
	\quad
	\begin{tikzpicture}[scale=0.3]
		\draw [thick](0,0)--(1,1)--(2,2)--(3,1)--(4,2)--(5,1);
		\draw[red,thick] (5,1)--(6,0);
	\end{tikzpicture}
	\quad
	\begin{tikzpicture}[scale=0.3]
		\draw [thick](0,0)--(1,1)--(2,0)--(3,1)--(4,0)--(5,1)--(6,0);
	\end{tikzpicture}
\end{equation*}
\caption{The 10 paths redrawn, with red south-east edges instead of south-west edges.}
	\end{figure}

The rules to generate such decorated Dyck paths are: each edge $(1,-1)$ may be  black or red, but
\begin{tikzpicture}[scale=0.3]\draw [thick](0,0)--(1,1); \draw [red,thick] (1,1)--(2,0);\end{tikzpicture}
and
\begin{tikzpicture}[scale=0.3] \draw [red,thick] (0,1)--(1,0);\draw [thick](1,0)--(2,1);\end{tikzpicture}
are forbidden.

Our interest is in particular in \emph{partial} decorated Dyck paths, ending at level $j$, for fixed $j\ge0$;
the instance $j=0$ is the classical case. 

The analysis of partial skew Dyck paths was recently  started in \cite{Baril-neu} (using the notion `prefix of a skew Dyck path') using
Riordan arrays instead of our kernel method. The latter  gives us \emph{bivariate} generating functions,
from which it is easier to draw conclusions. Two variables, $z$ and $u$, are used, where $z$ marks the length
of the path and $j$ marks the end-level. We briefly mention that one can, using a third variable $w$, also
count the number of red edges.

Again, once all generating functions are explicitly known, many corollories can be derived in a standard fashion.
We only do this in a few instances. But we would like to emphasize that the substitution
\begin{equation*}
x=\frac{v}{1+3v+v^2},
\end{equation*}
which was used in \cite{HPW, Prodinger-hex} allows to write \emph{explicit enumerations}, using
the notion of a (weighted) trinomial coefficient:
\begin{equation*}
\binom{n;1,3,1}{k}:=[t^k](1+3t+t^2)^n.
\end{equation*}

The second part of the paper deals with a dual version, where the paths are read from right to left.

\section{Generating functions and the kernel method}
\label{dunno}

We catch the essence of a decorated Dyck path using a state-diagram:

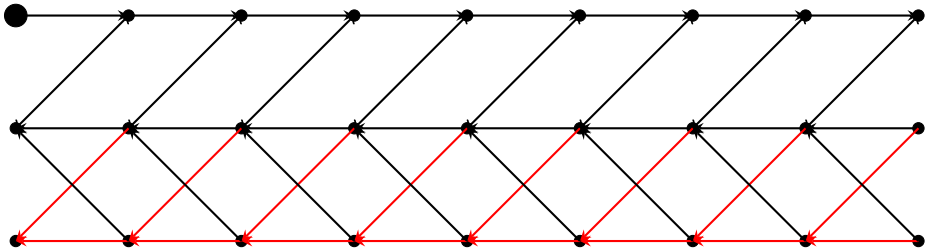
\begin{figure}[h]

\begin{center}
\begin{tikzpicture}[scale=1.5]
	\draw (0,0) circle (0.1cm);
	\fill (0,0) circle (0.1cm);
	
	\foreach \x in {0,1,2,3,4,5,6,7,8}
	{
		\draw (\x,0) circle (0.05cm);
		\fill (\x,0) circle (0.05cm);
	}

	\foreach \x in {0,1,2,3,4,5,6,7,8}
{
	\draw (\x,-1) circle (0.05cm);
	\fill (\x,-1) circle (0.05cm);
}

\foreach \x in {0,1,2,3,4,5,6,7,8}
{
	\draw (\x,-2) circle (0.05cm);
	\fill (\x,-2) circle (0.05cm);
}

\foreach \x in {0,1,2,3,4,5,6,7}
{
	\draw[ thick,-stealth] (\x,0) -- (\x+1,0);

}

\foreach \x in {0,1,2,3,4,5,6,7}
{
	\draw[ thick,<->] (\x+1,0) -- (\x,-1);
	
}

\foreach \x in {0,1,2,3,4,5,6,7}
{
	\draw[ thick,-stealth] (\x+1,-1) -- (\x,-1);
	
}
\foreach \x in {0,1,2,3,4,5,6,7}
{
	\draw[ thick,-stealth,red] (\x+1,-1) -- (\x,-2);
	
}

\foreach \x in {0,1,2,3,4,5,6,7}
{
	\draw[ thick,-stealth,red] (\x+1,-2) -- (\x,-2);
	
}

\foreach \x in {0,1,2,3,4,5,6,7}
{
	\draw[ thick,-stealth] (\x+1,-2) -- (\x,-1);
	
}

\end{tikzpicture}
\end{center}
\caption{Three layers of states according to the type of steps leading to them (up, down-black, down-red).}
\end{figure}
It has three types of states, with $j$ ranging from 0 to infinity; in the drawing, only $j=0..8$ is shown. The first
layer of states refers to an up-step leading to a state, the second layer refers to a black down-step leading to a state
and the third layer refers to a red down-step leading to a state. We will work out generating functions describing all paths
leading to a particular state. We will use the notations $f_j,g_j,h_j$ for the three respective layers, from top to bottom.
Note that the syntactic rules of forbidden patterns
\begin{tikzpicture}[scale=0.3]\draw [thick](0,0)--(1,1); \draw [red,thick] (1,1)--(2,0);\end{tikzpicture}
and
\begin{tikzpicture}[scale=0.3] \draw [red,thick] (0,1)--(1,0);\draw [thick](1,0)--(2,1);\end{tikzpicture}
can be clearly seen from the picture. The functions depend on the variable $z$ (marking the number of steps), but mostly
we just write $f_j$ instead of $f_j(z)$, etc.

The following recursions can be read off immediately from the diagram:
\begin{gather*}
f_0=1,\quad f_{i+1}=zf_i+zg_i,\quad i\ge0,\\
g_i=zf_{i+1}+zg_{i+1}+zh_{i+1},\quad i\ge0,\\
	h_i=zh_{i+1}+zg_{i+1},\quad i\ge0.
\end{gather*}
And now it is time to introduce the promised \emph{bivariate} generating functions: 
\begin{equation*}
F(z,u)=\sum_{i\ge0}f_i(z)u^i,\quad
G(z,u)=\sum_{i\ge0}g_i(z)u^i,\quad
H(z,u)=\sum_{i\ge0}h_i(z)u^i.
\end{equation*}
Again, often we just write $F(u)$ instead of $F(z,u)$ and treat $z$ as a `silent' variable. Summing the recursions leads to
\begin{align*}
\sum_{i\ge0}u^if_{i+1}&=\sum_{i\ge0}u^izf_i+\sum_{i\ge0}u^izg_i,\\
\sum_{i\ge0}u^ig_i&=\sum_{i\ge0}u^izf_{i+1}+\sum_{i\ge0}u^izg_{i+1}+\sum_{i\ge0}u^izh_{i+1},\\
\sum_{i\ge0}u^ih_i&=\sum_{i\ge0}u^izh_{i+1}+\sum_{i\ge0}u^izg_{i+1}.
\end{align*}
This can be rewritten as
\begin{align*}
	\frac1u(F(u)-1)&=zF(u)+zG(u),\\*
	G(u)&=\frac zu(F(u)-1)+\frac zu(G(u)-G(0))+\frac zu(H(u)-H(0)),\\*
	H(u)&=	\frac zu(G(u)-G(0))+\frac zu(H(u)-H(0)).
\end{align*}
This is a typical application of the kernel method. For a gentle example-driven introduction to the kernel method, see \cite{Prodinger-kernel}. First,
\begin{align*}
F(u)&=\frac{z^2uG(0)+z^2uH(0)+z^2u-u-z^3+2z}{-{z}^{3}-u+2z+z{u}^{2}-{z}^{2}u},\\
G(u)&=\frac{z(H(0)-uzH(0)+z^2+G(0)-zuG(0)-zu)}{-{z}^{3}-u+2z+z{u}^{2}-{z}^{2}u},\\
H(u)&=\frac{z(-uzH(0)-z^2-zuG(0)+G(0)-z^2H(0)+H(0)-z^2G(0))}{-{z}^{3}-u+2z+z{u}^{2}-{z}^{2}u}.
\end{align*}
The denominator factors as $z(u-r_1)(u-r_2)$, with
\begin{equation*}
r_1=\frac{1+z^2+\sqrt{1-6z^2+5z^4}}{2z},\quad r_2=\frac{1+z^2-\sqrt{1-6z^2+5z^4}}{2z}.
\end{equation*}
Note that $r_1r_2=2-z^2$. 
 Since the factor $u-r_2$ in the denominator is ``bad,'' it must also cancel in the numerators. From this
we conclude as a first step
\begin{equation*}
G(0) =  \frac{1-2z^2H(0)-3z^2-\sqrt{1-6z^2+5z^4}}{2z^2},
\end{equation*}
and by further simplification
\begin{equation*}
H(0)=\frac{1-4z^2+z^4+(z^2-1)\sqrt{1-6z^2+5z^4}}{2-z^2}.
\end{equation*}
Thus (with $W=\sqrt{1-6z^2+5z^4}=\sqrt{(1-z^2)(1-5z^2)}$\,)
\begin{align*}
F(u)&=\frac{-1-z^2-W}{2z(u-r_1)}=\frac{1+z^2+W}{2zr_1(1-u/r_1)},\\
G(u)&=\frac{-1+z^2+W}{2z(u-r_1)}=\frac{1-z^2-W}{2zr_1(1-u/r_1)},\\
H(u)&=\frac{-1+3z^2+W}{2z(u-r_1)}=\frac{1-3z^2-W}{2zr_1(1-u/r_1)}.
\end{align*}
The total generating function is
\begin{equation*}
S(u)=F(u)+G(u)+H(u)=\frac{3-3z^2-W}{2zr_1(1-u/r_1)}.
\end{equation*}
 The coefficient of $u^jz^n$ in $S(u)$ counts the partial paths of length $n$, ending at level $j$. 
 We will write $s_j=[u^j]S(u)$.
Furthermore
\begin{align*}
f_j=[u^j]	F(u)&=[u^j]\frac{1+z^2+W}{2zr_1(1-u/r_1)},\\
g_j=[u^j]	G(u)&=[u^j]\frac{1-z^2-W}{2zr_1(1-u/r_1)},\\
h_j=[u^j]	H(u)&=[u^j]\frac{1-3z^2-W}{2zr_1(1-u/r_1)}.
\end{align*}
At this stage, we are only interested in
\begin{equation*}
s_j=f_j+g_j+h_j=[u^j]\frac{3-3z^2-W}{2zr_1(1-u/r_1)}=\frac{3-3z^2-W}{2zr_1^{j+1}},
\end{equation*}
which is the generating function of all (partial) paths ending at level $j$. Parity considerations give us that only coefficients $[z^n]s_j$ are non-zero if $n\equiv j\bmod2$.
To make this more transparent, we set 
\begin{equation*}
P(z)=zr_1=\frac{1+z^2+\sqrt{1-6z^2+5z^4}}{2}, 
\end{equation*}
and then
\begin{equation*}
	s_j=f_j+g_j+h_j=z^j\frac{3-3z^2-W}{2P^{j+1}}.
\end{equation*}
Now we read off coefficients. We do this using residues and contour integration. The path of integration, in both variables $x$ resp.\ $v$ is a
small circle or an equivalent contour.
\begin{align*}
	[z^{2m+j}]s_j&=[z^{2m}]\frac{3-3z^2-W}{2P^{j+1}}=
	[x^m]\frac{3-3x-\sqrt{1-6x+5x^2}}{2\Big(\frac{1+x-\sqrt{1-6x+5x^2}}{2}\Big)^{j+1}}\\
	&=[x^m]\frac{3-3\frac v{1+3v+v^2}-\frac{1-v^2}{1+3v+v^2}}{2\big(\frac{v(v+2)}{1+3v+v^2}\big)^{j+1}}\\
	&=[x^m]\frac{(1+v)(1+2v)}{v^{j+1}(v+2)^{j+1}}(1+3v+v^2)^j\\
	&=\frac1{2\pi i}\oint\frac{dx}{x^{m+1}}\frac{(1+v)(1+2v)}{v^{j+1}(v+2)^{j+1}}(1+3v+v^2)^j\\
	&=\frac1{2\pi i}\oint\frac{dv}{v^{m+1}}\frac{(1+v)(1+2v)(1-v^2)}{v^{j+1}(v+2)^{j+1}}(1+3v+v^2)^{m-1+j}\\
	&=[v^{m+j+1}]\frac{(1+v)^2(1+2v)(1-v)}{(v+2)^{j+1}}(1+3v+v^2)^{m-1+j}.
\end{align*}
Note that
\begin{equation*}(1+v)^2(1+2v)(1-v)=
-9+27( v+2 ) -29( v+2 ) ^{2}+13( v+2) ^{3}-2( v+2 ) ^{4};
\end{equation*}
consequently
\begin{align*}
[v^k]&\frac{(1+v)^2(1+2v)(1-v)}{(v+2)^{j+1}}\\
&=-9\frac1{2^{j+1+k}}\binom{-j-1}{k}
+27\frac1{2^{j+k}}\binom{-j}{k}
-29\frac1{2^{j-1+k}}\binom{-j+1}{k}\\&
+13\frac1{2^{j-2+k}}\binom{-j+2}{k}
-2\frac1{2^{j-3+k}}\binom{-j+3}{k}=:\lambda_{j;k}.
\end{align*}
With this abbreviation we find
\begin{equation*}
	[v^{m+j+1}]\frac{(1+v)^2(1+2v)(1-v)}{(v+2)^{j+1}}(1+3v+v^2)^{m-1+j}
	=\sum_{k=0}^{m+j+1}\lambda_{j;k}\binom{m-1+j;1,3,1}{m+j+1-k}.
\end{equation*}
This is not extremely pretty but it is \emph{explicit} and as good as it gets.
Here are the first few generating functions:
\begin{align*}
s_0&=1+z^2+3z^4+10z^6+36z^8+137z^{10}+543z^{12}+\cdots\\*
s_1&=z+2z^3+6z^5+21z^7+79z^9+311z^{11}+1265z^{13}+\cdots\\
s_2&=z^2+3z^4+10z^6+37z^8+145z^{10}+589z^{12}+2455z^{14}+\cdots\\
s_3&=z^3+4z^5+15z^7+59z^9+241z^{11}+1010z^{13}+4314^{15}+\cdots\\
\end{align*}
We could also give such lists for the functions $f_j$, $g_j$, $h_j$, if desired. We summarize the essential findings of this section:
\begin{theorem} The generating function of decorated (partial) Dyck paths, consisting of $n$ steps, ending on level $j$, is given by
	\begin{equation*}
S(z,u)=\frac{3-3z^2-\sqrt{1-6z^2+5z^4}}{2zr_1(1-u/r_1)},
	\end{equation*}
	 with
	 \begin{equation*}
r_1=\frac{1+z^2+\sqrt{1-6z^2+5z^4}}{2z}.
	 \end{equation*}
 Furthermore
 \begin{equation*}
[u^j]S(z,u)=\frac{3-3z^2-\sqrt{1-6z^2+5z^4}}{2zr_1^{j+1}}.
 \end{equation*}
\end{theorem}

\section{Open ended paths}

If we do not specify the end of the paths, in other words we sum over all $j\ge0$, then at the level of generating functions
this is very easy, since we only have to set $u:=1$.
We find
\begin{align*}
S(1)&=-\frac{(z+1)(z^2+3z-2)+(z+2)\sqrt{1-6z^2+5z^4}}{2z(z^2+2z-1)}\\
&=1+z+2z^2+3z^3+7z^4+11z^5+26z^6+43z^7+102z^8+175z^9+416z^{10}+\cdots.
\end{align*}

\section{Counting red edges}

We can use an extra variable, $w$, to count additionally the red edges that occur in a path. We use the same
letters for generating functions. Eventually, the coefficient $[z^nu^jw^k]S$ is the number of (partial) paths consisting of $n$ steps, leading
to level $j$, and having passed $k$ red edges. The endpoint of the original skew path has then coordinates $(n-2k,j)$. The computations are very similar, and we only
sketch the key steps.

\begin{equation*}
	f_0=1,\quad f_{i+1}=zf_i+zg_i,\quad i\ge0,
\end{equation*}
\begin{equation*}
	g_i=zf_{i+1}+zg_{i+1}+zh_{i+1},\quad i\ge0,
\end{equation*}
\begin{equation*}
	h_i=wzh_{i+1}+wzg_{i+1},\quad i\ge0;
\end{equation*}
\begin{align*}
	\frac1u(F(u)-1)&=zF(u)+zG(u),\\*
	G(u)&=\frac zu(F(u)-1)+\frac zu(G(u)-G(0))+\frac zu(H(u)-H(0)),\\*
	H(u)&=	\frac {wz}u(G(u)-G(0))+\frac {wz}u(H(u)-G(0));
\end{align*}
\begin{align*}
	F(u)&=\frac{z^2uG(0)+z^2uH(0)+z^2u-u-wz^3+z+wz}{-w{z}^{3}-u+z+wz+z{u}^{2}-w{z}^{2}u},\\
	G(u)&=\frac{z(H(0)-uzH(0)+wz^2+G(0)-zuG(0)-zu)}{-w{z}^{3}-u+z+wz+z{u}^{2}-w{z}^{2}u},\\
	H(u)&=\frac{wz(-uzH(0)-z^2-zuG(0)+G(0)-z^2H(0)+H(0)-z^2G(0))}{-w{z}^{3}-u+z+wz+z{u}^{2}-w{z}^{2}u}.
\end{align*}
The denominator factors as $z(u-r_1)(u-r_2)$, with
\begin{align*} 
	r_1&=\frac{1+wz^2+\sqrt{1-(4+2w)z^2+(4w+w^2)z^4}}{2z},\\*
	r_2&=\frac{1+wz^2-\sqrt{1-(4+2w)z^2+(4w+w^2)z^4}}{2z}.
\end{align*}
Note the factorization $1-(4+2w)z^2+(4w+w^2)z^4=(1-z^2w)(1-(4+w)z^2)$. 
Since the factor $u-r_2$ in the denominator is ``bad,'' it must also cancel in the numerators. From this
we eventually find, with the abbreviation
$W=\sqrt{1-(4+2w)z^2+(4w+w^2)z^4}\,$)
\begin{align*}
	F(u)&=\frac{-1-wz^2-W}{2z(u-r_1)},\\
	G(u)&=\frac{-1+wz^2+W}{2z(u-r_1)},\\
	H(u)&=\frac{-1+(2+w)z^2+W}{2z(u-r_1)}.
\end{align*}
The total generating function is
\begin{equation*}
	S(u)=F(u)+G(u)+H(u)=\frac{-2-w+z^2(w+w^2)+ wW}{2z(u-r_1)}.
	\end{equation*}
The special case $u=0$ (return to the $x$-axis) is to be noted:
\begin{equation*}
	S(0)=\frac{-2-w+z^2(w+w^2)+ wW}{-2zr_1}=\frac{1-wz^2-W}{2z^2}.
\end{equation*}
Since there are only even powers of $z$ in this function, we replace $x=z^2$ and get
\begin{align*}
	S(0)&=\frac{1-wx-\sqrt{1-(4+2w)x+(4w+w^2)x^2}}{2x}\\
	&=1+x+(w+2)x^2+(w^2+4w+5)x^3+(w^3+6w^2+15w+14)x^4+\cdots.
\end{align*}
Compare the factor $(w^2+4w+5)$ with the earlier drawing of the 10 paths.

There is again a substitution that allows for better results:
\begin{equation*}
z=\frac{v}{1+(2+w)v+v^2}, \quad\text{then}\quad S(0)=1+v.
\end{equation*}
Reading off coefficients can now be done using modified trinomial coefficients:
\begin{equation*}
\binom{n;1,2+w,1}{k}=[t^k]\bigl(1+(2+w)t+t^2\bigr)^n.
\end{equation*}
Again, we use contour integration to extract coefficients:
\begin{align*}
[x^n](1+v)&=\frac1{2\pi i}\oint \frac{dx}{x^{n+1}}(1+v)\\
&=\frac1{2\pi i}\oint \frac{dx}{v^{n+1}}\frac{1-v^2}{(1+(2+w)v+v^2)^2}(1+(2+w)v+v^2)^{n+1}(1+v)\\
&=[v^n](1-v)(1+v)^2(1+(2+w)v+v^2)^{n-1}\\
&=\binom{n-1;1,2+w,1}{n}+\binom{n-1;1,2+w,1}{n-1}\\*
&\qquad-\binom{n-1;1,2+w,1}{n-2}-\binom{n-1;1,2+w,1}{n-3}.
\end{align*}

Now we want to count the average number of red edges. For that, we differentiate $S(0)$ w.r.t.\ $w$, followed by $w:=1$.
This leads to
\begin{equation*}
\frac{-1+6x-5x^2+(1+3x)\sqrt{1-6x+5x^2}}{2(1-x)(1-5x)}.
\end{equation*}

A simple application of singularity analysis leads to
\begin{equation*}
\frac{\frac1{2\sqrt5}[x^n]\frac1{\sqrt{1-5x}}}{-\sqrt5[x^n]\sqrt{1-5x}}\sim \frac {n}{5}.
\end{equation*}
So, a random path consisting of $2n$ steps has about $n/5$ red steps, on average. 

For readers who are not familiar with singularity analysis of generating functions \cite{FlOd90, FS}, we just
mention that one determines the local expansion around the dominating singularity, which is at $z=\frac15$ in our instance.
In the denominator, we just have the total number of skew Dyck paths, according to the sequence A002212 in \cite{OEIS}.

In the example of Figure~2, the exact average is $6/10$, which curiously is exactly the same as $3/5$.

We finish the discussion by considering fixed powers of $w$ in $S(0)$, counting skew Dyck paths consisting of zero, one, two, three, \dots red edges. We find
\begin{align*}
[w^0]S(0)&=\frac{1-\sqrt{1-4x}}{2x},\\
[w^1]S(0)&=\frac{1-2x-\sqrt{1-4x}}{2\sqrt{1-4x}},\\
[w^2]S(0)&=\frac{x^3}{(1-4x)^{3/2}},\\
[w^3]S(0)&=\frac{x^4(1-2x)}{(1-4x)^{5/2}},\\
[w^4]S(0)&=\frac{x^5(1-4x+5x^2)}{(1-4x)^{7/2}}, \quad\&\text{c}.
\end{align*}
The generating function $[w^0]S(0)$ is of course the generating function of Catalan numbers, since no red edges just means: ordinary Dyck paths.
We can also conclude that the asymptotic behaviour is of the form $n^{k-3/2}4^n$, where the polynomial contribution gets higher, but the exponential growth
stays the same: $4^n$. This is compared to the scenario of an \emph{arbitrary} number of red edges, when we get an exponential growth of the form $5^n$.

\section{Dual skew Dyck paths}

The mirrored version of skew Dyck paths with two types of up-steps, $(1,1)$ and $(-1,1)$ are also cited among the objects in A002212 in \cite{OEIS}.
We call them dual skew paths and drop the `dual' when it isn't necessary. When the paths come back to the $x$-axis, no new enumeration is necessary, but this
is no longer true for paths ending at level $j$.

Here is a list of the 10 skew paths consisting of 6 steps:

\begin{figure}[h]
	\begin{equation*}
		\begin{tikzpicture}[scale=0.3]
			\draw [thick](0,0)--(3,3)--(6,0);
		\end{tikzpicture}
		\quad
		\begin{tikzpicture}[scale=0.3]
			\draw [thick](0,0)--(-1,1)--(1,3)--(4,0);
		\end{tikzpicture}
		\quad
		\begin{tikzpicture}[scale=0.3]
			\draw [thick](0,0)--(1,1)--(0,2)--(1,3)--(4,0);
		\end{tikzpicture}
		\quad
		\begin{tikzpicture}[scale=0.3]
			\draw [thick](0,0)--(-2,2)--(-1,3)--(2,0);
		\end{tikzpicture}
		\quad
		\begin{tikzpicture}[scale=0.3]
			\draw [thick](0,0)--(2,2)--(4,0)--(5,1)--(6,0);
		\end{tikzpicture}
	\end{equation*}
	\begin{equation*}
		\begin{tikzpicture}[scale=0.3]
			\draw [thick](0,0)--(1,1)--(2,0)--(4,2)--(6,0);
		\end{tikzpicture}
		\quad
		\begin{tikzpicture}[scale=0.3]
			\draw [thick](0,0)--(-1,1)--(0,2)--(2,0)--(3,1)--(4,0);
		\end{tikzpicture}
		\quad
		\begin{tikzpicture}[scale=0.3]
			\draw [thick](0,0)--(1,1)--(2,2)--(3,1)--(4,2)--(6,0);
		\end{tikzpicture}
		\quad
		\begin{tikzpicture}[scale=0.3]
			\draw [thick](0,0)--(-1,1)--(0,2)--(1,1)--(2,2)--(4,0);
		\end{tikzpicture}
		\quad
		\begin{tikzpicture}[scale=0.3]
			\draw [thick](0,0)--(1,1)--(2,0)--(3,1)--(4,0)--(5,1)--(6,0);
		\end{tikzpicture}
	\end{equation*}
	\caption{All 10 dual skew Dyck paths of length 6 (consisting of 6 steps).}
\end{figure}

We prefer to work with the equivalent model (resembling more traditional Dyck paths) where
we replace each step $(-1,-1)$ by $(1,-1)$ but label it blue. Here is the list of the 10 paths again (Figure 2):

\begin{figure}[h]
	\begin{equation*}
		\begin{tikzpicture}[scale=0.3]
			\draw [thick](0,0)--(3,3)--(6,0);
		\end{tikzpicture}
		\quad
		\begin{tikzpicture}[scale=0.3]
			\draw [thick,cyan](0,0)--(1,1);
			\draw [thick](1,1)--(3,3)--(6,0);
		\end{tikzpicture}
		\quad
		\begin{tikzpicture}[scale=0.3]
			\draw [thick](0,0)--(1,1);
			\draw [thick,cyan](1,1)--(2,2);
			\draw [thick](2,2)--(3,3)--(6,0);
		\end{tikzpicture}
		\quad
		\begin{tikzpicture}[scale=0.3]
			\draw [thick,cyan](0,0)--(2,2);
			\draw [thick](2,2)--(3,3)--(6,0);
		\end{tikzpicture}
		\quad
		\begin{tikzpicture}[scale=0.3]
			\draw [thick](0,0)--(2,2)--(4,0)--(5,1)--(6,0);
		\end{tikzpicture}
	\end{equation*}
	\begin{equation*}
		\begin{tikzpicture}[scale=0.3]
			\draw [thick](0,0)--(1,1)--(2,0)--(4,2)--(6,0);
		\end{tikzpicture}
		\quad
		\begin{tikzpicture}[scale=0.3]
			\draw [thick,cyan](0,0)--(1,1);
			\draw [thick](1,1)--(2,2)--(4,0)--(5,1)--(6,0);
		\end{tikzpicture}
		\quad
		\begin{tikzpicture}[scale=0.3]
			\draw [thick](0,0)--(1,1)--(2,2)--(3,1)--(4,2)--(6,0);
		\end{tikzpicture}
		\quad
		\begin{tikzpicture}[scale=0.3]
			\draw [thick,cyan](0,0)--(1,1);
			\draw [thick](1,1)--(2,2)--(3,1)--(4,2)--(6,0);
		\end{tikzpicture}
		\quad
		\begin{tikzpicture}[scale=0.3]
			\draw [thick](0,0)--(1,1)--(2,0)--(3,1)--(4,0)--(5,1)--(6,0);
		\end{tikzpicture}
	\end{equation*}
	\caption{All 10 dual skew Dyck paths of length 6 (consisting of 6 steps).}
\end{figure}

The rules to generate such decorated Dyck paths are: each edge $(1,-1)$ may be  black or blue, but
\begin{tikzpicture}[scale=0.3]\draw [thick](0,1)--(1,0); \draw [cyan,thick] (1,0)--(2,1);\end{tikzpicture}
and
\begin{tikzpicture}[scale=0.3] \draw [cyan,thick] (0,0)--(1,1);\draw [thick](1,1)--(2,0);\end{tikzpicture}
are forbidden.

Our interest is in particular in \emph{partial} decorated Dyck paths, ending at level $j$, for fixed $j\ge0$;
the instance $j=0$ is the classical case. 

The analysis of partial skew Dyck paths was recently  started in \cite{Baril-neu} (using the notion `prefix of a skew Dyck path') using
Riordan arrays instead of our kernel method. The latter  gives us \emph{bivariate} generating functions,
from which it is easier to draw conclusions. Two variables, $z$ and $u$, are used, where $z$ marks the length
of the path and $j$ marks the end-level. We briefly mention that one can, using a third variable $w$, also
count the number of blue edges.

The substitution
\begin{equation*}
	x=\frac{v}{1+3v+v^2},
\end{equation*}
which was used in \cite{HPW, Prodinger-hex} is the key to the success and allows to write \emph{explicit enumerations}, using
the notion of a (weighted) trinomial coefficient:
\begin{equation*}
	\binom{n;1,3,1}{k}:=[t^k](1+3t+t^2)^n.
\end{equation*}

\section{Generating functions and the kernel method}

We catch the essence of a decorated (dual skew) Dyck path using a state-diagram:

\begin{figure}[h]
	
	\begin{center}
		\begin{tikzpicture}[scale=1.5]
			\draw (0,0) circle (0.1cm);
			\fill (0,0) circle (0.1cm);
			
			\foreach \x in {0,1,2,3,4,5,6,7,8}
			{
				\draw (\x,0) circle (0.05cm);
				\fill (\x,0) circle (0.05cm);
			}
			
			\foreach \x in {0,1,2,3,4,5,6,7,8}
			{
				\draw (\x,-1) circle (0.05cm);
				\fill (\x,-1) circle (0.05cm);
			}
			
			\foreach \x in {0,1,2,3,4,5,6,7,8}
			{
				\draw (\x,1) circle (0.05cm);
				\fill (\x,1) circle (0.05cm);
			}
			
			\foreach \x in {0,1,2,3,4,5,6,7}
			{
				\draw[ thick,-latex] (\x,0) -- (\x+1,0);
				
			}
			\foreach \x in {0,1,2,3,4,5,6,7}
			{
				\draw[ thick,-latex,cyan] (\x,1) -- (\x+1,1);
				
			}
			
			\foreach \x in {0,1,2,3,4,5,6,7}
			{
				\draw[ thick,-latex,cyan] (\x,0) to (\x+1,1);
				
			}
			
			\foreach \x in {0,1,2,3,4,5,6,7}
			{
				\draw[ thick,-latex] (\x,1)  to (\x+1,0);
				
			}
			
			\foreach \x in {0,1,2,3,4,5,6,7}
			{
				\draw[ thick,latex-latex] (\x+1,0)  to (\x,-1);
				
			}
			
			\foreach \x in {0,1,2,3,4,5,6,7}
			{
				\draw[ thick,-latex] (\x+1,-1)  to (\x,-1);
				
			}

		\end{tikzpicture}
	\end{center}
	\caption{Three layers of states according to the type of steps leading to them (down, up-black, up-blue).}
\end{figure}
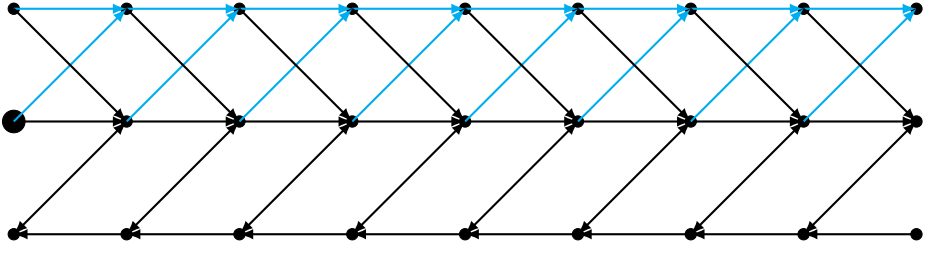
It has three types of states, with $j$ ranging from 0 to infinity; in the drawing, only $j=0..8$ is shown. The first
layer of states refers to an up-step leading to a state, the second layer refers to a black down-step leading to a state
and the third layer refers to a blue down-step leading to a state. We will work out generating functions describing all paths
leading to a particular state. We will use the notations $c_j,a_j,b_j$ for the three respective layers, from top to bottom.
Note that the syntactic rules of forbidden patterns
\begin{tikzpicture}[scale=0.3]\draw [thick,cyan](0,0)--(1,1); \draw [thick] (1,1)--(2,0);\end{tikzpicture}
and
\begin{tikzpicture}[scale=0.3] \draw [thick] (0,1)--(1,0);\draw [thick,cyan](1,0)--(2,1);\end{tikzpicture}
can be clearly seen from the picture. The functions depend on the variable $z$ (marking the number of steps), but mostly
we just write $a_j$ instead of $a_j(z)$, etc.

The following recursions can be read off immediately from the diagram:
\begin{gather*}
	a_0=1,\quad a_{i+1}=za_i+zb_i+zc_i,\quad i\ge0,\\
	b_i=za_{i+1}+zb_{i+1},\quad i\ge0,\\
	c_{i+1}=za_{i}+zc_{i},\quad i\ge0.
\end{gather*}
And now it is time to introduce the promised \emph{bivariate} generating functions: 
\begin{equation*}
	A(z,u)=\sum_{i\ge0}a_i(z)u^i,\quad
	B(z,u)=\sum_{i\ge0}b_i(z)u^i,\quad
	C(z,u)=\sum_{i\ge0}c_i(z)u^i.
\end{equation*}
Again, often we just write $A(u)$ instead of $A(z,u)$ and treat $z$ as a `silent' variable. Summing the recursions leads to
\begin{align*}
	\sum_{i\ge0}u^ia_i &=1+u\sum_{i\ge0}u^i(za_i+zb_i+zc_i)\\
	&=1+uzA(u)+uzB(u)+uzC(u),\\
	\sum_{i\ge0}u^ib_i &= \sum_{i\ge0}u^i(za_{i+1}+zb_{i+1})\\
	&=\frac zu\sum_{i\ge1}u^ia_i+\frac zu\sum_{i\ge1}u^ib_i,\\
	\sum_{i\ge1}u^ic_i &=uz\sum_{i\ge0}u^ia_i+uz\sum_{i\ge0}u^ic_i.
\end{align*}
This can be rewritten as
\begin{align*}
	A(u)&=1+uzA(u)+uzB(u)+uzC(u),\\
	B(u)&=\frac zu(A(u)-a_0)+\frac zu(B(u)-b_0),\\
	C(u)&=c_0+uzA(u)+uzC(u).
\end{align*}
Note that $a_0=1$, $c_0=0$.
Simplification leads to
\begin{equation*}
	C(u)=\frac{uzA(u)}{1-uz}
\end{equation*}
and
\begin{equation*}
	B(u)=\frac{z(A(u)-1-B(0))}{u-z}
\end{equation*}
leaving us with just one equation
\begin{equation*}
	A(u)={\frac { \left( z-u+u{z}^{2}+u{z}^{2}B(0) \right)  \left( uz-1
			\right) }{{u}^{2}{z}^{3}+u{z}^{2}-2{u}^{2}z-z+u}}.
\end{equation*}
This is a typical application of the kernel method, \cite{Prodinger-kernel}. 
\begin{equation*}
	{u}^{2}{z}^{3}+u{z}^{2}-2{u}^{2}z-z+u=z(z^2-2)(u-s_1)(u-s_2)
\end{equation*}
The denominator factors as $2z(z^2-2)(u-s_1)(u-s_2)$, with
\begin{equation*}
	s_1=\frac{1+z^2+\sqrt{1-6z^2+5z^4}}{2z(2-z^2)},\quad s_2=\frac{1+z^2-\sqrt{1-6z^2+5z^4}}{2z(2-z^2)}.
\end{equation*}
Note that $s_1s_2=\frac{1}{2-z^2}$. 
Since the factor $u-s_2$ in the denominator is ``bad,'' it must also cancel in the numerators. From this
we conclude (again with the abbreviation $W=\sqrt{1-6z^2+5z^4}\,$)
\begin{equation*}
	B(0) =  \frac{zs_2}{1-2zs_2},
\end{equation*}
and further
\begin{equation*}
	A(u)
	=\frac{(1-uz)(1+z^2+W)}{2z(z^2-2)(u-s_1)},
\end{equation*}
\begin{equation*}
	B(u)=\frac{1-2z^2-W}{z(2-z^2)(u-s_1)},
\end{equation*}
\begin{equation*}
	C(u)=\frac{1+z^2+W}{2(z^2-2)}\frac{u}{u-s_1},
\end{equation*}
and for the function of main interest
\begin{equation*}
	G(u)=A(u)+B(u)+C(u)=\frac{3z^2-3+W}{2z(2-z^2)(u-s_1)}.
\end{equation*}
Note that
\begin{align*}
	\frac1{s_1}&=\frac{1+z^2-\sqrt{1-6z^2+5z^4}}{2z}=zS,\\
	\frac1{s_2}&=\frac{1+z^2+\sqrt{1-6z^2+5z^4}}{2z}.
\end{align*}
Then
\begin{align*}
	[u^j]G(u)&=[u^j]\frac{3z^2-3+W}{2z(z^2-2)s_1(1-u/s_1)}\\
	&=\frac{3z^2-3+W}{2z(z^2-2)s_1^{j+1}}
	=\frac{3z^2-3+W}{2(z^2-2)}z^{j}S^{j+1}.
\end{align*}
So $[u^j]G(u)$ contains only powers of the form $z^{j+2N}$. Now we continue
\begin{align*}
	[z^{j+2N}u^j]G(u)&
	=[z^{2N}]\frac{3z^2-3+W}{2(z^2-2)}S^{j+1}
	\\&=[x^{N}]\frac{3x-3+\sqrt{1-6x+5x^2}}{2(x-2)}\bigg(\frac{1+x-\sqrt{1-6x+5x^2}}
	{2x}\bigg)^{j+1}\\
	&=[x^{N}](v+1)(v+2)^{j}
\end{align*}
which is the generating function of all (partial) paths ending at level $j$.

Now we read off coefficients. We do this using residues and contour integration. The path of integration, in both variables $x$ resp.\ $v$ is a
small circle or an equivalent contour;
\begin{align*}
	[z^{j+2N}u^j]G(u)&=[x^{N}](v+1)(v+2)^{j}\\
	&=\frac1{2\pi i}\oint \frac{dx}{x^{N+1}}(v+1)(v+2)^{j}\\
	&=\frac1{2\pi i}\oint \frac{dv}{v^{N+1}}(1+3v+v^2)^{N+1}\frac{(1-v^2)}{(1+3v+v^2)^2}(v+1)(v+2)^{j}\\
	&=[v^{N}](1+3v+v^2)^{N-1}(1-v)(1+v)^2(v+2)^{j}.
\end{align*}
Note that
\begin{equation*}(1-v)(1+v)^2=
	3-7( v+2 ) +5( v+2 ) ^{2}- ( v+2) ^{3};
\end{equation*}
consequently
\begin{align*}
	[z^{j+2N}u^j]G(u)&=[v^{N}](1+3v+v^2)^{N-1}\Big[3-7( v+2 ) +5( v+2 ) ^{2}- ( v+2) ^{3}
	\Big](v+2)^{j}.
\end{align*}
We abbreviate:
\begin{align*}
	\mu_{j;k}&=[v^{k}]\Big[3(v+2)^{j}-7(v+2)^{j+1} +5(v+2)^{j+2}- (v+2)^{j+3}\Big]\\
	&=3\binom{j}{k}2^{j-k}-7\binom{j+1}{k}2^{j+1-k}+5\binom{j+2}{k}2^{j+2-k}-\binom{j+3}{k}2^{j+3-k}.
\end{align*}
With this notation we get
\begin{equation*}
	[z^{j+2N}u^j]G(u)
	=\sum_{0\le k\le N-1}\mu_{j;k}\binom{N-1;1,3,1}{N-k}.
\end{equation*}
Here are the first few generating functions:
\begin{align*}
	G_0&=1+{z}^{2}+3{z}^{4}+10{z}^{6}+36{z}^{8}+137{z}^{10}+543{z}^{
		12}+2219{z}^{14}
	+\cdots\\*
	G_1&=2z+3{z}^{3}+10{z}^{5}+36{z}^{7}+137{z}^{9}+543{z}^{11}+
	2219{z}^{13}+9285{z}^{15}
	+\cdots\\
	G_2&=4{z}^{2}+8{z}^{4}+29{z}^{6}+111{z}^{8}+442{z}^{10}+1813{z
	}^{12}+7609{z}^{14}+32521{z}^{16}
	+\cdots\\
	G_3&=8{z}^{3}+20{z}^{5}+78{z}^{7}+315{z}^{9}+1306{z}^{11}+5527
	{z}^{13}+23779{z}^{15}+103699{z}^{17}
	+\cdots\\
\end{align*}
We could also give such lists for the functions $a_j$, $b_j$, $c_j$, if desired. We summarize the essential findings of this section:
\begin{theorem} The generating function of decorated (partial) dual skew Dyck paths, consisting of $n$ steps, ending on level $j$, is given by
	\begin{equation*}
		G(z,u)=\frac{3z^2-3+\sqrt{1-6z^2+5z^4}}{2z(2-z^2)(u-s_1)},
	\end{equation*}
	with
	\begin{equation*}
		s_1=\frac{2z}{1+z^2-\sqrt{1-6z^2+5z^4}}.
	\end{equation*}
	Furthermore
	\begin{equation*}
		[u^j]G(z,u)=\frac{3z^2-3+\sqrt{1-6z^2+5z^4}}{2(z^2-2)}z^jS^{j+1},
	\end{equation*}
	with
	\begin{equation*}
		S=\frac{1+z^2-\sqrt{1-6z^2+5z^4}}{2z^2}.
	\end{equation*}
\end{theorem}

\section{Open ended paths}

If we do not specify the end of the paths, in other words we sum over all $j\ge0$, then at the level of generating functions
this is very easy, since we only have to set $u:=1$.
We find
\begin{align*}
	G(1)&=\frac{(1+z)(1-3z)}{2z(z^2+2z-1)-\sqrt{1-6z^2+5z^4}}\\
	&=1+2z+5{z}^{2}+11{z}^{3}+27{z}^{4}+62{z}^{5}+151{z}^{6}+
	354{z}^{7}+859{z}^{8}+2036{z}^{9}+\cdots.
\end{align*}

\section{Counting blue edges}

We can use an extra variable, $w$, to count additionally the blue edges that occur in a path. We use the same
letters for generating functions. Eventually, the coefficient $[z^nu^jw^k]S$ is the number of (partial) paths consisting of $n$ steps, leading
to level $j$, and having passed $k$ blue edges. The endpoint of the original skew path has then coordinates $(n-2k,j)$. The computations are very similar, and we only
sketch the key steps.

\begin{gather*}
	a_0=1,\quad a_{i+1}=za_i+zb_i+zc_i,\quad i\ge0,\\
	b_i=za_{i+1}+zb_{i+1},\quad i\ge0,\\
	c_{i+1}=wza_{i}+wzc_{i},\quad i\ge0.
\end{gather*}
This leads to
\begin{align*}
	A(u)&=1+uzA(u)+uzB(u)+uzC(u),\\
	B(u)&=\frac zu(A(u)-a_0)+\frac zu(B(u)-b_0),\\
	C(u)&=c_0+wuzA(u)+wuzC(u).
\end{align*}
Solving,
\begin{equation*}
	S(u)=A(u)+B(u)+C(u)={\frac {u-wu{z}^{2}-zA(0)-zB(0)+uw{z}^{2}A(0)+uw{z}^{2}B(0)}{{u}^{2}{z}^{3}w+u-w{u}^{2}z-{u}^{2}z-z+wu{z}^{2}}}.
\end{equation*}
The denominator factors as $-z(1+w-z^2w)(u-s_1)(u-s_2)$, with
\begin{align*} 
	s_1&={\frac {1+{z}^{2}w+\sqrt {1-2\,{z}^{2}w+{z}^{4}{w}^{2}-4\,{z}^{2
				}+4{z}^{4}w}}{2z \left( 1+w-{z}^{2}w \right) }},\\*
	s_2&= {\frac {1+{z}^{2}w-\sqrt {1-2\,{z}^{2}w+{z}^{4}{w}^{2}-4\,{z}^{2
				}+4{z}^{4}w}}{2z \left( 1+w-{z}^{2}w \right) }}	.
\end{align*}
Note the factorization $1-(4+2w)z^2+(4w+w^2)z^4=(1-z^2w)(1-(4+w)z^2)$. 
Since the factor $u-r_2$ in the denominator is ``bad,'' it must also cancel in the numerators. From this
we eventually find, with the abbreviation
$W=\sqrt{1-(4+2w)z^2+(4w+w^2)z^4}\,$)
\begin{equation*}
	G(0)={\frac {1-{z}^{2}w-W }{2{z}^{2}}},
\end{equation*}
and further
\begin{equation*}
	G(u)=\frac {w-{z}^{2}{w}^{2}-wW+2-2{z}^{2}w}
	{2z \left( -w	-1+{z}^{2}w \right) (u-s_1)}.
\end{equation*}
The special case $u=0$ (return to the $x$-axis) is to be noted:
\begin{equation*}
	G(0)=1+{z}^{2}+ \left( w+2 \right) {z}^{4}+ \left( {w}^{2}+4w+5 \right) 
	{z}^{6}+ \left( w+2 \right)  \left( {w}^{2}+4w+7 \right) {z}^{8}+\cdots.
\end{equation*}
Compare the factor $(w^2+4w+5)$ with the earlier drawing of the 10 paths.
There is again a substitution that allows for better results:
\begin{equation*}
	z=\frac{v}{1+(2+w)v+v^2}, \quad\text{then}\quad G(0)=1+v.
\end{equation*}
Since $S(u)=G(u)$ with $S(u)$ from the first part of the paper, as it means the same objects, read from left to right resp.\ from right to left, no 
new analysis is required.

\section{Skew paths that can go into negative territory}

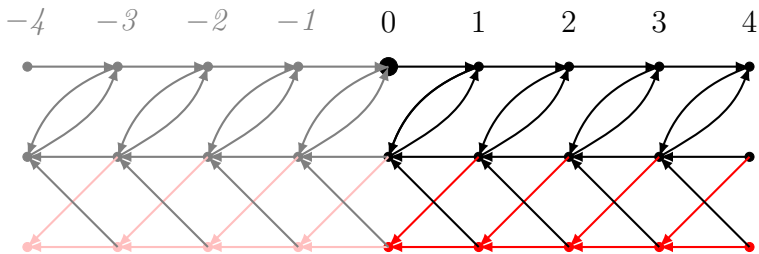
\begin{figure}[h]

	\begin{center}
		\begin{tikzpicture}[scale=1.2]
			\draw (0,0) circle (0.1cm);
			\fill (0,0) circle (0.1cm);
			
			\foreach \x in {0,1,2,3,4}
			{
				\draw (\x,0) circle (0.05cm);
				\fill (\x,0) circle (0.05cm);
			}
			\foreach \x in {1,2,3,4}
			{
				\draw[gray] (-\x,0) circle (0.05cm);
				\fill[gray] (-\x,0) circle (0.05cm);
			}

			\foreach \x in {0,1,2,3,4}
			{
				\draw (\x,-1) circle (0.05cm);
				\fill (\x,-1) circle (0.05cm);
			}
			
			\foreach \x in {1,2,3,4}
			{
				\draw [gray](-\x,-1) circle (0.05cm);
				\fill[gray] (-\x,-1) circle (0.05cm);
			}
			
			\foreach \x in {1,2,3,4}
			{
				\draw [pink](-\x,-2) circle (0.05cm);
				\fill[pink] (-\x,-2) circle (0.05cm);
			}
			
			\foreach \x in {0,1,2,3,4}
			{
				\draw[red] (\x,-2) circle (0.05cm);
				\fill [red](\x,-2) circle (0.05cm);
			}
			
			\foreach \x in {0,1,2,3}
			{
				\draw[ thick,-latex] (\x,0) -- (\x+1,0);
				
			}
			
			\foreach \x in {0,1,2,3}
			{
				\draw[ thick,-latex,gray] (\x-4,0) -- (\x+1-4,0);
				
			}
			
			\foreach \x in {1,2,3}
			{
				\draw[thick,  -latex] (\x+1,0) to[out=200,in=70]  (\x,-1);	
			}
			
			\foreach \x in {0,1,2,3}
			{
				
				\draw[thick,  -latex,gray] (\x+1-4,0) to[out=200,in=70]  (\x-4,-1);	
			}

			\draw[thick,  -latex] (1,0) to[out=200,in=70]  (0,-1);

			\foreach \x in {0,1,2,3}
			{
				
				\draw[thick,  -latex] (\x,-1) to[out=30,in=250]  (\x+1,0);	
				
			}
			\draw[thick,  -latex] (1,0) to[out=200,in=70]  (0,-1);

			\foreach \x in {0,1,2,3}
			{
				
				\draw[thick,  -latex,gray] (\x-4,-1) to[out=30,in=250]  (\x+1-4,0);	
				
			}
			
			\foreach \x in {0,1,2,3}
			{
				\draw[ thick,-latex] (\x+1,-1) -- (\x,-1);
				
			}
			\foreach \x in {0,1,2,3}
			{
				\draw[ thick,-latex,gray] (\x+1-4,-1) -- (\x-4,-1);
				
			}
			
			\foreach \x in {0,1,2,3}
			{
				\draw[ thick,-latex,red] (\x+1,-1) -- (\x,-2);
				
			}
			
			\foreach \x in {0,1,2,3}
			{
				\draw[ thick,-latex,pink] (\x+1-4,-1) -- (\x-4,-2);
				
			}
			
			\foreach \x in {0,1,2,3}
			{
				\draw[ thick,-latex,red] (\x+1,-2) -- (\x,-2);
				
			}
			
			\foreach \x in {0,1,2,3}
			{
				\draw[ thick,-latex,pink] (\x+1-4,-2) -- (\x-4,-2);
				
			}
			
			\foreach \x in {0,1,2,3}
			{
				\draw[ thick,-latex] (\x+1,-2) -- (\x,-1);
				
			}
			
			\foreach \x in {0,1,2,3}
			{
				\draw[ thick,-latex,gray] (\x+1-4,-2) -- (\x-4,-1);
				
			}
			
			\foreach \x in {0,1,2,3,4}
			{
				\draw  (\x,0.5) node  {$\x$};

			}
			\foreach \x in {-1,-2,-3,-4}
			{
				\draw  (\x,0.5) node[gray]  {$\mathit \x$};
				
			}

		\end{tikzpicture}
	\end{center}
	\caption{Three layers of states according to the type of steps leading to them (up, down-black, down-red).}
	\label{arseneg}
\end{figure}

For Dyck paths and the standard random walk on the integers, the enumeration, if the negative territory is allowed, is easier.
In our instance of paths equipped with an additional red down-step and the usual restrictions (up--red and red--up are forbidden) this is
not so; it is rather more complicated. The paths may be described by another directed graph~Figure~\ref{arseneg}.

We have the following recursions,
\begin{gather*}
	f_{i}=[i=0]+zf_{i-1}+zg_{i-1},\\
	g_i=zf_{i+1}+zg_{i+1}+zh_{i+1},\\
	h_i=zg_{i+1}+zh_{i+1}.
\end{gather*}
For negative indices we need to introduce separate sequences,
\begin{equation*}
	a_i=f_{-i}, \ b_i=g_{-i}, \ c_i=h_{-i}.
\end{equation*}
Then we find
\begin{gather*}
	f_{-i}=[-i=0]+zf_{-i-1}+zg_{-i-1},\\*
	g_{-i}=zf_{-i+1}+zg_{-i+1}+zh_{-i+1},\\*
	h_{-i}=zg_{-i+1}+zh_{-i+1}
\end{gather*}
and, rewriting,
\begin{gather*}
	a_{i}=[i=0]+za_{i+1}+zb_{i+1},\\
	b_{i}=za_{i-1}+zb_{i-1}+zc_{i-1},\\
	c_{i}=zb_{i-1}+zc_{i-1}.
\end{gather*}
Introducing
\begin{equation*}
F(u)=\sum_{i\ge0}f_iu^i,\ G(u)=\sum_{i\ge0}g_iu^i,\ H(u)=\sum_{i\ge0}h_iu^i
\end{equation*}
and
\begin{equation*}
	A(u)=\sum_{i\ge0}a_iu^i,\ B(u)=\sum_{i\ge0}b_iu^i,\ C(u)=\sum_{i\ge0}c_iu^i
\end{equation*}
we get the following 6 equations:
\begin{align*}
	F(u)-f_0&=zu(F(u)+G(u)),\\
	G(u)&=\frac zu(F(u)+G(u)+H(u)-f_0-g_0-h_0),\\
	H(u)&=\frac zu(G(u)+H(u)-g_0-h_0),\\
	A(u)&=1+\frac zu(A(u)+B(u)-f_0-g_0),\\
	B(u)-g_0&=zu(A(u)+B(u)+C(u)),\\
	C(u)-h_0&=zu(B(u)+C(u)).
\end{align*}
Solving the system,
\begin{align*}
	F&={\frac {{z}^{2}ug_0+{z}^{2}uh_0+{z}^{2}u{f_0}-u{f_0}-{z}^{3}{f_0}+2 z{f_0}}{-{z}^{3}-u+2 z+z{u}^{2}-{z}^{2}u}},\\
	G&=-{\frac {z \left( zu{h_0}-{z}^{2}{f_0}-{g_0}-{h_0}
			+uz{f_0}+zu{g_0} \right) }{-{z}^{3}-u+2 z+z{u}^{2}-{z}^{2}u}},\\
	H&=-{\frac {z \left( zu{h_0}+{z}^{2}{f_0}+zu{g_0}+{z}^{2}{
				h_0}-{g_0}+{z}^{2}{g_0}-{h_0} \right) }{-{z}^{3}-u
			+2 z+z{u}^{2}-{z}^{2}u}},\\
	A&={\frac {2 {z}^{2}u{f_0}+{z}^{2}u{g_0}+{z}^{2}u{h_0}-2
			z{u}^{2}-z{f_0}+u}{{z}^{2}u+{z}^{3}{u}^{2}+u-2 z{u}^{2}-z}},\\
	B&={\frac {z{u}^{2}-{z}^{2}u{f_0}+u{g_0}-z{g_0}-{z}^{2}{u}
			^{3}+{z}^{3}{u}^{2}{f_0}+{z}^{3}{u}^{2}{g_0}-z{u}^{2}{g_0}+z{u}^{2}{h_0}-{z}^{2}u{h_0}}{{z}^{2}u+{z}^{3}{u}^{2}+u
			-2 z{u}^{2}-z}},\\
	C&=-{\frac {-{z}^{2}{u}^{3}+{z}^{3}{u}^{2}{f_0}+{z}^{3}{u}^{2}{g_0}-u{h_0}-z{u}^{2}{g_0}+z{u}^{2}{h_0}+z{h_0}+{z
			}^{2}u{g_0}}{{z}^{2}u+{z}^{3}{u}^{2}+u-2 z{u}^{2}-z}},
\end{align*}
and
\begin{equation*}
	-{z}^{3}-u+2 z+z{u}^{2}-{z}^{2}u=z(u-r_1)(u-r_2)
\end{equation*}
with
\begin{equation*}
	r_{1,2}=\frac{1+z^2\pm\sqrt{1-6z^2+5z^4}}{2z}.
\end{equation*}
As usual, the factor $u-r_2$ must cancel out. The other denominators are
\begin{equation*}
{z}^{2}u+{z}^{3}{u}^{2}+u
-2 z{u}^{2}-z=z(z^2-2)(u-s_1)(u-2_2)
\end{equation*}
and $s_1=1/r_2$, $s_2=1/r_1$. The factor $u-s_1$ must cancel out as well. This leads to
\begin{align*}
	F(u)&=\frac {{z}^{2}g_0+{z}^{2}h_0+{z}^{2}f_0-f_0}{r_2z-1-{z}^{2}+zu},\\
	G(u)&=-{\frac { \left(f_0+ g_0+h_0 \right) {z}^{2}}{r_2z-1-{z}^{2}+zu}},\\
	H(u)&=-{\frac {{z}^{2} \left( g_0+h_0 \right) }{r_2z-1-{z}^{2}+zu}},
\end{align*}
\scriptsize
\begin{align*}
	A(u)&=-{\frac {2 s_1z-2 {z}^{2}f_0-{z}^{2}g_0-{z}^{2}h_0-1+2 zu}{s_1{z}^{3}-2s_1z+1+{z}^{2}+u{z}^{3}-2 zu}},\\
	B(u)&=-\frac{\mathcal{X}}{s_1 {z}^{3}
			-2 s_1 z+1+{z}^{2}+u{z}^{3}-2 zu},\\
	C(u)&={\frac {{s_1}^{2}{z}^{2}-s_1 {z}^{3}f_0-s_1 {z}^{3}g_0+s_1 zg_0-s_1 zh_0+s_1 {z}^{2}u+h-{z}^{2}g_0-u{z}^{3}f_0
			-u{z}^			{3}g_0+zug_0-zuh_0+{z}^{2}{u}^{2}}{s_1 {z}^{3}-2 s_1 z+1+{z}^{2			}+u{z}^{3}-2 zu}}
\end{align*}
\normalsize with $\mathcal{X}= {s_1}^{2}{z}^{2}-s_1 z-s_1 {z}^{3}f_0-s_1
	{z}^{3}g_0+s_1 zg_0-s_1 zh_0+s_1 {z}^{2}u+{z}^{2}f_0-g_0+{z}
	^{2}h_0-zu-u{z}^{3}f_0-u{z}^{3}g_0+zug_0-zuh_0+{z}^{2}{u}^{2}$.
The computation of $f_0$, $g_0$, $h_0$ requires some care.
From the equations for $G$ and $H$ we conclude
\begin{equation*}
	g_0=\frac{z^2f_0+h_0}{1-z^2}
\end{equation*}
and from the expression for $H$, as just derived, we find
\begin{equation*}
	h_0=f_0{\frac {{z}^{4}}{r_2{z}^{3}-r_2z+1-2{z}^{2}}}.
\end{equation*}
So both, $g_0$ and $h_0$ are multiples of $f_0$. As any $f_0$ would solve the first 3 equations with the appropriate $g_0$, $h_0$,
we need to resort to $A(u)$ since there we find that  $f_0=1+\cdots$. By elimination,
\begin{equation*}
	A(u)={\frac {-2z{u}^{2}+u+g_0{z}^{2}u+{z}^{2}uh_0+2f_0{z}^{2}u-zf_)}{{z}^{2}u-z+{z}^{3}{u}^{2}-2z{u}^{2}+u}}.
\end{equation*}
Substitute $u=0$ and use $g_0$ and $h_0$ from before and then solve $A(0)=f_0$ leads to
\begin{equation*}
	f_0={\frac {r_2{z}^{3}-r_2z+1-2{z}^{2}}{-{z}^{6}+{z}^{4}+2r_2{z}^{3}-3{z}^{2}-r_2z+1}}.
\end{equation*}
It can be made explicit:
\begin{equation*}
f_0=\frac{1+z^2-\sqrt{1-6z^2+5z^4}}{2z^2(2-z^2)}.
\end{equation*}
Now everything is explicit:
\begin{align*}
	f_0&=1+z^2+2z^4+6z^6+21z^8+79z^{10}+311z^{12}+1265z^{14},\\
	g_0&=z^2+3z^4+10z^6+37z^8+145z^{10}+589z^{12}+2455z^{14},\\
	h_0&=z^4+5z^6+21z^8+87z^{10}+365z^{12}+1555z^{14}.
\end{align*}
The expressions for $g_0$ and $h_0$ are a bit long, but
\begin{align*}
	f_0+g_0+h_0&=\frac{1-3z^2+2z^4-\sqrt{1-6z^2+5z^4}}{2z^4(2-z^2)}\\
	&=1+2z^2+6z^4+21z^6+79z^8+311z^{10}+1265z^{12}+5275z^{14}.
\end{align*}
The coefficients  $1,2,6,21,\dots$ are sequence A033321 in \cite{OEIS}. In the comments to this sequence, 
the number of skew Dyck paths of semilength $n$ ending with a down step $(1,-1)$ is mentioned, something that follows from
our results for $g_0$ in section~\ref{dunno}.


\bibliographystyle{plain}


\end{document}